\documentclass[a4paper,12pt]{article}
\pdfoutput=1

\usepackage[margin=2cm]{geometry}
\usepackage[english]{babel}
\usepackage[latin1]{inputenc}
\usepackage{times}
\usepackage[T1]{fontenc}
\usepackage{graphicx}
\usepackage{xcolor}
\usepackage{amsmath}
\usepackage{amssymb}
\usepackage{amsbsy}
\usepackage{trfsigns}
\usepackage[round]{natbib}
\usepackage{tikz}

\newtheorem{assumption}{Assumption}
\newtheorem{theorem}{Theorem}

\newtheorem{definition}{Definition}

\def\j{\mathbf{j}}

\def\0{\mathbf{0}}

\def\sgn{\mathrm{sgn}}
\def\diag{\mathrm{diag}}

\title{
Generalization of the construction method\footnote{Parts of the results here have been presented in a preliminary version, Multistability equivalence between gene regulatory networks of different dimensionality, which appeared in the \textit{Proceedings of the 12th European Control Conference (ECC)}; see~\citet{Schittler2013b}.}\\
for multistability-equivalent gene regulatory networks\\
to systems with multi-input multi-output loopbreaking
}
\author{Daniella Schittler$^1$, Taouba Jouini$^1$, Frank Allg{\"o}wer$^1$, Steffen Waldherr$^2$ \\
\small{$^1$Institute for Systems Theory and Automatic Control, University of Stuttgart}\\  %$^a$ 
\small{$^2$Institute for Automation Engineering, Otto-von-Guericke University Magdeburg}\\
\small{daniella.schittler@ist.uni-stuttgart.de, taoubaJouini@gmail.com,} \\
\small{steffen.waldherr@ovgu.de, frank.allgower@ist.uni-stuttgart.de}
}

\begin{document}

\maketitle

%%%%%%%%%%%%%%%
%%%%%%%%%%%%%%%
%%%%%%%%%%%%%%%
%%%%%%%%%%%%%%%

\begin{abstract}
The problem of equivalence in terms of multistability properties between gene regulatory network models of different dimensionality has been recently addressed by~\citet{Schittler2013b}. The authors in that work proposed construction rules for a high-dimensional dynamical system, when given a low-dimensional dynamical system and the high-dimensional network structure.
However, the proof therein was restricted to the class of systems for which all internal feedback loops can be broken by a loopbreaking approach yielding a single-input single-output (SISO) system.
In this report, we present the generalization of the proof to systems with any number of internal feedback loops, which will be broken by a generalized loopbreaking approach resulting in a multi-input multi-output (MIMO) system.
This generalization of the method renders the construction method applicable to a broad class of gene regulatory network models, thus promoting the transfer of results from core motif models to more realistic, high-dimensional models of gene regulation.
We demonstrate the potential and value of our method by applying it to an example of a gene regulatory network in mesenchymal stem cell differentiation.
\end{abstract}

%%%%%%%%%%%%%%%
%%%%%%%%%%%%%%%
\section{Introduction}
\label{sec:ch3:intro}

The dynamics of gene regulatory networks (GRNs) are commonly studied via nonlinear ordinary differential equation (ODE) models.
Generic modeling approaches make use of the finding that a small number of ``master regulator'' genes and ``core motifs'' are largely responsible for generating multistable or oscillatory behavior~(\citealp{Alon2007,Huang2007,Tyson2003}).
Corresponding core motif models of low dimensionality can be developed and analyzed based on systems-theoretic tools such as multistability and bifurcation analysis (see, for example,~\citet{Huang2005,Schittler2010}).
However, such low-dimensional core motif models are meant for the conceptual understanding of higher-dimensional gene regulatory systems. 
As typically tens or hundreds of genes may be involved in the determination of a cell type, the aim is to derive also suitable high-dimensional GRN models that are closer to reality and can be fit to experimental data measuring the expression time course from multiple genes.

Recently, a construction method was proposed~(\citealp{Schittler2013b}) that contributes to closing the gap between low-dimensional and high-dimensional GRN models: 
First, the authors introduced the concept of multistability equivalence which allows to match the stability properties of two GRN systems of different dimensionality.
Second, the construction method proposed therein was proven to yield a high-dimensional GRN system that is multistability-equivalent to a given low-dimensional GRN system, while consistent with a given interaction network structure.
However, this result was restricted by the assumption that all feedback loops in the system can be broken by a single-input single-output (SISO) loopbreaking, meaning that a loopbreaking~(\citealp{Waldherr2009}) exists that results in a SISO-system which will not have any unstable pole-zero cancellations any more.

In this report, we generalize the respective step in the proof to GRN systems with arbitrary many feedback loops in the interaction function. 
We present a loopbreaking that guarantees to break all internal dynamics except for the degradation self-loops, which results in a multi-input multi-output (MIMO) system.
This loopbreaking approach is used to conduct the proof and thus show multistability equivalence for a general class of GRN systems.
The MIMO loopbreaking approach as well as the generalized proof was first sketched in~\citet{Jouini2013}, whereas the aim of this report is to put the proof step in direct relation to the original proof of~\citet{Schittler2013b}, and to make the generalized proof broadly available.

In Section~\ref{sec:problem}, we briefly recall the concept of multistability equivalance.
Then we review in Section~\ref{sec:construction} the construction method of~\citet{Schittler2013b}.
In Section~\ref{sec:proof}, we recall the theorem stated in~\citet{Schittler2013b}, and then present in detail the generalized step in the proof regarding the loopbreaking.
%Section~\ref{sec:bifu} makes a statement on bifurcations that can be directly drawn from the multistability equivalence.
In Section~\ref{sec:example}, we apply the proposed method to a GRN in mesenchymal stem cell differentiation.
Finally, Section~\ref{sec:summ} concludes with a short summary and outlook.

%%%%%%%%%%%%%%%
\section{Definition of multistability equivalence} % and problem formulation
\label{sec:problem}

In this section we review the definition of~\citet{Schittler2013b} for a multistability-equivalent system as well as assumptions imposed by~\citet{Schittler2013b}, and recall the problem addressed therein.

To capture the problem, a low-dimensional system is defined for which the dynamics (and thus also its structure) are given, and a high-dimensional system is defined for which only the structure, but importantly no dynamics, are known.

Let be given an ODE system 
\begin{equation}
\Sigma_f: \dot{z} = f(z) = a(z) - d(z), \ z \in \mathbb{R}_{+}^n,
\label{eq:sys_f}
\end{equation}
with the interaction rate $a(z) \in \mathbb{R}^n_+$, 
and the degradation rate $d(z) = k z,\ k=\diag(k_1,\ldots,k_n)$, with $k_i > 0$.

We assume the system to have $R$ steady states $z^{\ast (r)}, r=1\ldots R$ which are determined by $f(z^{\ast (r)})=0$.
The system's Jacobian at a steady state $z^{\ast (r)}$ is denoted by 
\begin{align}
J_f^{(r)} := \frac{\partial f}{\partial z}(z^{\ast (r)})
.
\label{eq:jacobian}
\end{align}
 
The structure of the low-dimensional system will be represented by an interaction sign matrix,
\begin{align}
S_a := \sgn\left( \frac{\partial a}{\partial z} \right) \in \{-1,0,+1\}^{n \times n}
.
\label{eq:S_a}
\end{align}

The number of unstable modes (positive eigenvalues) for each steady state is given by
\begin{align}
|\{ \lambda^{(r)}_{f,v} | \Re(\lambda^{(r)}_{f,v}) > 0 \} |  ,
\label{eq:lambda_f}
\end{align}
with $\lambda^{(r)}_{f,v}, v=1\ldots n$ the eigenvalues of the Jacobian at the steady state,  $J_f|_{z=z^{\ast (r)}}$. 

Let an interaction sign matrix of a higher-dimensional system be given as
\begin{equation}
S_A \in \{-1,0,+1\}^{N \times N}, N \geq n.
\label{eq:S_A}
\end{equation}

The following definition is made in order to formulate the concept of a 
higher-dimensional ODE system 
\begin{equation}
\Sigma_F: \dot{x} = F(x) = A(x) - D(x), \ x \in \mathbb{R}_{+}^N, %x \geq 0,
\label{eq:sys_F}
\end{equation}
which is desired to have equivalent multistability properties.
\begin{definition}
A system~\eqref{eq:sys_F}
is called an \textit{$N$-dimensional multistability-equivalent} system to \eqref{eq:sys_f} and \textit{consistent with the interaction structure} given by \eqref{eq:S_A}, if the following hold.
\begin{itemize}
	\item[(i)] The derivative of the interaction function $A(x)$ has signs as given by the interaction sign matrix \eqref{eq:S_A}:
	\begin{equation}
	\sgn\left( \frac{\partial A}{\partial x} \right) = S_A \ \forall x \in \mathbb{R}_{++}^N.
	\label{eq:sgnA}
	\end{equation}
	 %$\sgn( \frac{\partial A}{\partial x}) = S_A$ as given by \eqref{eq:S_A}.
	 %
	 %
	\item[(ii)] There exists an injective map 
	%$h: z^{\ast} \mapsto x^{\ast}=h(z^{\ast}) \in \mathbb{R}^N$ with:
	$h: \mathbb{R}^n_+ \rightarrow \mathbb{R}^N_+, z^{\ast} \mapsto h(z^{\ast})$, with:
	\vspace{-0.35cm}
	\begin{equation}
	f(z^{\ast}) = 0 \Leftrightarrow F(h(z^{\ast})) = 0.
	\label{eq:equiv_ii}
	\end{equation}
	\item[(iii)] The number of unstable modes (positive eigenvalues) in both systems \eqref{eq:sys_f} and \eqref{eq:sys_F}, for each pair of steady states $z^{\ast (r)}$, $x^{\ast (r)} = h(z^{\ast (r)})$, is equal:
	\begin{equation}
	|\{ \lambda^{(r)}_{F,u} | \Re(\lambda^{(r)}_{F,u}) > 0 \} | = 
	|\{ \lambda^{(r)}_{f,v} | \Re(\lambda^{(r)}_{f,v}) > 0 \} |.
	\label{eq:equiv_iii}
	\end{equation}
\end{itemize}
\label{def:ms-eq}
\end{definition}
In the remainder, a system that meets the properties of Def.~\ref{def:ms-eq} will be called \textit{multistability-equivalent} for short.

Further technical assumptions on the systems are imposed as follows:
For the Jacobians $J_f^{(r)}$, $r=1,\ldots,R$ as in~\eqref{eq:jacobian} it is assumed that they do not have eigenvalues on the imaginary axis. 
This is important for the structural stability of the system, especially when the argument principle will be exploited later in the proof.

The interaction sign matrix of the low-dimensional system~\eqref{eq:S_a}
is assumed to be constant over $z \in \mathbb{R}_{++}^n$, 
since otherwise it will depend on the specific values of the state variables $z$ and thus the system structure is not uniquely determined.

The interaction sign matrix of the higher-dimensional system~\eqref{eq:S_A} may be available, for example, from qualitative knowledge about gene interactions. These can be represented by an interaction graph, which in turn can be directly translated into a sign matrix~\eqref{eq:S_A}.

The class of considered systems is restricted for technical reasons, by the two following assumtions:
\begin{assumption}
\textit{(Modular structure)}
Assume that, possibly by reordering the state space variables of system \eqref{eq:sys_F}, the interaction sign matrix \eqref{eq:S_A} fulfills the following structural property.
There exist numbers $m_i \in \{1, \dotsc, n\}$, $i=n+1,\dotsc,N$, such that 
\begin{itemize}
	\item in rows $i = n+1,\ldots,N$, columns $j = 1,\ldots,n$: $S_A|_{i,m_i} = +1$ for no more than one $m_i$.
	\item in rows $i = n+1,\ldots,N$, columns $j = n+1,\ldots,N$: $S_A|_{i,j} \in \{0,+1\}$ for all $j$ where $m_i = m_j$, and $S_A|_{i,j} = 0$ otherwise.
\end{itemize}
\label{ass:0.0}
\end{assumption}

In the remainder, the first $n$ state variables $x_i,\ i \in \{1,\ldots,n\}$ will be referred to as ``master genes'', whereas the remaining $(N-n)$ state variables $\{x_{n+1},\ldots,x_N\}$ will be referred to as ``module genes''. 
In this way, each index $i \in \{n+1,\ldots,N\}$ is uniquely assigned to one interaction module: $i \in \mathcal{M}_{m_i}$, 
being disjoint subsets of the indices $\{n+1,\ldots,N\} \supseteq \mathcal{M}_j: \ \bigcap_j \mathcal{M}_j = \emptyset$, 
and with the assignment such that 
\begin{itemize}
	\item interactions affecting $x_i$ come either from $x_{m_i}$, or from other $x_j$ belonging to the same interaction module,
	\item interactions from $x_i$ go to genes that are master genes or belong to the same module, 
	\item interactions between genes belonging to the same module are nonnegative. 
\end{itemize}
The general structure of such a matrix is given in Figure~\ref{fig:matrixstructure}.

\begin{figure}%[!ht]
%\scalebox{0.8}{
\begin{small}
\begin{align*}
\begin{split}  
S_A &= 
\begin{matrix}
i=  & 								\\
%% rows 1...n
 	\begin{matrix}
	1 \\ \vdots \\ n
	\end{matrix}
	&
	\begin{tabular}{|c|}
	\hline \\
	%\in \{-1,0,+1\}^{n \times N} \\
	%\{-1,0,+1\}^{n \times N} \\
	\hspace{2.55cm} -1,0,+1  \hspace{2.55cm} \\
	\\
	\hline
	\end{tabular}
	\\
	%\hline
%% rows n+1...N
 	\begin{matrix}
	n+1 \\ \vdots \\ \vdots \\ \vdots \\ \vdots \\ \vdots \\ N
	\end{matrix}
	&
 %% = all interactions affecting module genes
 \begin{matrix}
	% the interactions from the module master gene {n+1,..,N} x {1,..,n}
	\begin{tabular}{|c|c|c|}
	\hline 
	  & &  \\
		$P_1$	 & 0 & 0  \\
	  & &  \\
	\hline 
	  & &  \\
	0 & $P_{\ldots}$ & 0 \\
	  & &  \\
	\hline 
	  & &  \\
	0 & 0 & $P_n$ \\
	  & &  \\
	\hline
	\end{tabular}
	&
	% the interactions within modules {n+1,..,N} x {n+1,..,N}
	\begin{tabular}{|c|c|c|}
	\hline 
	  & &  \\
	 \ 0,+1 \ & 0 & 0  \\
	  & &  \\
	\hline 
	  & &  \\
	0 & \ 0,+1 \ & 0  \\
	  & &  \\
	\hline 
	  & &  \\
	0 & 0 & \ 0,+1 \ \\
	  & &  \\
	\hline
	\end{tabular}
 \end{matrix}
\end{matrix}
\end{split}
\end{align*}
\end{small}
%}
\caption{General structure of a matrix fulfilling the structural requirements.
Each $P_m, m=1\ldots n$ is a a column vector of 0 and +1.}
\label{fig:matrixstructure}
\end{figure}

The classification into master and module genes with an according network structure 
might be predetermined from biological knowledge. 
If not, the selection of master genes could be addressed as a separate problem which is beyond the scope of this contribution.

\begin{assumption}
\textit{(Consistency of sign matrices $S_a$ and $S_A$)}
For each $(i,j) \in \{1,\ldots,n\} \times \{1,\ldots,n\}$ with $S_a|_{i,j} \neq 0$, there must exist some simple path from $j \in \{n+1,..,N\}$ to $i \in \{1,..,n\}$, %(a sequence of edges from j to i)
\begin{align}
\begin{split}    % j->i
p_{ij} &:= ( (i,\iota_{1}),(\iota_{1},\iota_{2}),\ldots,
(\iota_{\omega_{ij}-1},\iota_{\omega_{ij}}),(\iota_{\omega_{ij}},j) ) \\
& \quad \subseteq (\{1,\ldots,N\} \times \{1,\ldots,N\})^{(\omega_{ij}+1)}, \\ %\text{  for which } \\
& \text{ for which } \;  \prod_{(\iota',\iota'') \in p_{ij}} S_A|_{\iota',\iota''}  = S_a|_{i,j}.
\end{split}
%\label{eq:step0.1.a}
\end{align}
Moreover, if there exists $\iota_1 \in \{n+1,\ldots,N\}$ on some path from $j$ to $i$: $(i,\iota_1) \in p_{ij} $, %p_{ij_{\iota_1}}}
with $S_A|_{i,\iota_1} \neq 0$, there has to exist an interaction 
in the low-dimensional system
%$S_a|_{i,j} \neq 0$. %$: 
\begin{equation}
S_a|_{i,j} \neq 0. 
%\label{eq:step0.1.b}
\end{equation}
\label{ass:0.1}
\end{assumption}

If one of these assumptions is not fulfilled, then the proposed construction of an $N$-dimensional multistability-equivalent system can not be done in the way as proposed here, and therefore the existence of an $N$-dimensional multistability-equivalent system cannot be guaranteed.

%%%%%%%%%%%%%%%
\section{Construction rules}
\label{sec:construction}

Next, we review the construction method presented in~\citet{Schittler2013b}.
The aim of this construction is to solve the following problem:
Given a dynamic gene switch model~\eqref{eq:sys_f}, 
and the structure of a gene regulatory network, via an interaction sign matrix~\eqref{eq:S_A}, 
construct a dynamic gene regulatory network model~\eqref{eq:sys_F}, 
such that the obtained system~\eqref{eq:sys_F} is multistability-equivalent to the given system~\eqref{eq:sys_f} and consistent with the interaction structure~\eqref{eq:S_A}.

The idea of the construction is as follows: 
Additional interactions %(for the $N-n$ additional state variables) 
are introduced via linear activation functions.
The remaining interactions are defined in terms of the interactions from the low-dimensional system, with specific mappings between the state spaces of different dimensionality. 
These interactions are constructed such that the steady state gains of the additional interactions are exactly compensated when the system is at steady state.

The proposed construction procedure is as follows.
%% THE ACTUAL CONSTRUCTION:

%%% 1
\paragraph{Step 1:}
Construct functions $F_i(x)$ for the module genes indices $i = n+1,\ldots,N$, as follows: %for all $i \in \{n+1,\ldots,N\}$:
\begin{align}
\begin{split}
F_i(x) &= A_i(x) - D_i(x) \quad \text{with}\\
A_i(x) &= \sum_{j=1\ldots N} S_A|_{i,j} x_j,\\
D_i(x) &= K_i x_i, \; \text{with}\ K_i \in \mathbb{R}_{++} \; \text{to be chosen}. %determined}.
\end{split}
\label{eq:step1}
\end{align}
With this, interaction functions for the $(N-n)$ module genes are determined up to the parameters $K_i, i\in \{n+1,\ldots,N\}$.

%%% 2
\paragraph{Step 2:} Next, the influence of master genes on the module genes is captured by defining the following transfer gains. 
For all $(i,k),\ k \in \{1,\ldots,n\},\ i \in \mathcal{M}_k \subseteq \{n+1,\ldots,N\}$,
given a system \eqref{eq:sys_F}, denote the transfer gain with input $x_k$ and output $x_i$, via the system matrix $(\partial F / \partial x|_{l,m})_{l\in\mathcal{M}_k,m\in\mathcal{M}_k}$, that is, restricted to the state variables in the module $\mathcal{M}_k$:
%\begin{small}
\begin{align}
\begin{split}
G_{k \rightarrow i}&(\lambda) =
(e_i^T|_m)_{m\in\mathcal{M}_k} 
\left( \lambda I_{|\mathcal{M}_k|} - \left(\frac{\partial F}{\partial x}|_{l,m}\right)_{\substack{l\in\mathcal{M}_k\\m\in\mathcal{M}_k}} \right)^{-1}
(S_A|_{l,k})_{l\in\mathcal{M}_k} .
\end{split}
\label{eq:Gki}
\end{align}
%\end{small}
%
Then define parameters $\gamma_{ik}$ which give the corresponding steady state gain for input $x_k$ and output $x_i$:
\begin{align}
\begin{split}
\gamma_{ik} &:= G_{k \rightarrow i}(0)  \\
&=
(e_i^T|m)_{m\in\mathcal{M}_k} %\cdot
\left( - \left(\frac{\partial F}{\partial x}|_{l,m}\right)_{\substack{l\in\mathcal{M}_k\\m\in\mathcal{M}_k}} \right)^{-1} %\cdot
(S_A|_{l,k})_{l\in\mathcal{M}_k} 
.
\end{split}
\label{eq:gamma}
\end{align}
With this step, the signal transmission through the additionally introduced module genes is captured.

%%% 3
\paragraph{Step 3:} 
In this step, functions $F_i(x)$ for the master genes indices $i = 1,\ldots,n$ are constructed. 

For the signal transmission between master genes, it is desired to capture which interactions in the high-dimensional system have corresponding interactions in the low-dimensional system with either same or opposite sign. 
Therefore, for all pairs of genes, $(i,\nu)$ for $i,\nu  \in \{ 1,\ldots,n \}$, all interactions in the high-dimensional system from the master gene $\nu$ and its module $\mathcal{M}_{\nu}$ to the master gene $i$ are screened for having 
the same or opposite sign as the direct interaction from $\nu$ to $i$ in the low-dimensional system.

The results are captured in the $2n^2$ $N$-dimensional index vectors $J_{(i,\nu)}^{=}$ and $J_{(i,\nu)}^{\neq}$ that reflect the interactions having the same, respectively of opposite, sign when comparing the high-dimensional to the low-dimensional system's interaction structure. 
These index vectors $J_{(i,\nu)}^{=}$ and $J_{(i,\nu)}^{\neq}$, 
$i,\nu = 1,\ldots,n$, are defined as follows:
The $j$-th entry of $J_{(i,\nu)}^=$ is 
\begin{equation}
  \label{eq:J=}
  J_{(i,\nu)}^{=}|_{j} := 
  \left\{
  \begin{array}{l l}
  0  &, \text{ if } j \notin \mathcal{M}_{\nu} \cup\{\nu\}    \\
  \big|S_a|_{i,\nu} \big| \cdot \delta(S_A|_{i,j}, S_a|_{i,\nu})  &,
  	\text{ if } j \in \mathcal{M}_{\nu} \cup\{\nu\}
  \end{array}
  \right. ,
\end{equation}
with $\delta(a,b)$ the Kronecker delta of $a$ and $b$. 
Similarly, the $j$-th entry of $J_{(i,\nu)}^{\neq}$ is 
\begin{equation}
  \label{eq:Jneq}
  J_{(i,\nu)}^{\neq}|_{j} := 
  \left\{
  \begin{array}{l l}
  0  &, \text{ if } j \notin \mathcal{M}_{\nu} \cup\{\nu\}   \\
  \big|S_a|_{i,\nu} \big| \cdot \delta(-S_A|_{i,j}, S_a|_{i,\nu})  &,
  	\text{ if } j \in \mathcal{M}_{\nu} \cup\{\nu\}
  \end{array}
  \right. ,
\end{equation}
These vectors $J_{(i,\nu)}^{=}$ ($J_{(i,\nu)}^{\neq}$)
have an entry $1$ in the component $j$ whenever the interaction of $z_{\nu}$ onto $z_i$ in the low-dimensional system has the same (opposite, respectively) sign as the interaction of $x_j$ on $x_i$ in the high-dimensional system, whereas $x_j$ is within the module $\mathcal{M}_{\nu}$ of $x_{\nu}$, and zero otherwise.

Furthermore, we define an auxiliary map $\mu_i: \mathbb{R}^N \rightarrow \mathbb{R}^n$ to map from the higher-dimensional state space to the lower-dimensional state space. 
For this purpose, let define parameters $\gamma_{ik}$ similarly to \eqref{eq:gamma} but for the interactions between master genes, $i,k \in \{1,\ldots,n\}$:
\begin{align}
\begin{split}
\gamma_{ik} &:= \big| S_A|_{i,k} \big| . 
\end{split}
\label{eq:gamma_master}
\end{align}

Then, the auxiliary map $\mu_i: \mathbb{R}^N \rightarrow \mathbb{R}^n$ is defined as follows:
\begin{align}
\begin{split}
 & x \mapsto \mu_i(x) := 
 \Big( 
   \left( \frac{1 + \epsilon_{i}||J_{(i,\nu)}^{\neq}||}{||J_{(i,\nu)}^{=}||} \right)
      (J_{(i,\nu)}^{=})^T  (\gamma_{j\nu}^{-1} x_j)_{j=1\ldots N} 
      %\\
 %& \hspace{2cm} 
 - \frac{\epsilon_{i}}{||J_{(i,\nu)}^{=}||}(J_{(i,\nu)}^{\neq})^T  (\gamma_{j\nu}^{-1} x_j)_{j=1\ldots N}
 \Big)_{\nu=1\ldots n}
\end{split}
\label{eq:step3_mui}
\end{align}
%(or $\epsilon_{ij}$, but then rewrite into vector...) \\
with $\gamma_{j\nu} \in \mathbb{R}_{++}$ as determined in \eqref{eq:gamma}, \eqref{eq:gamma_master}.
Now, let the functions $F_i(x)$, $i = 1,\dotsc,n$ be constructed as follows, using the preceding definitions of index vectors and auxiliary map:
\begin{align}
\begin{split}
F_i(x) &= A_i(x) - D_i(x) \quad \text{with}\\
A_i(x) &= a_i(\mu_{i}(x)),\; \text{and}\\
D_i(x) &= k_i x_i. 
\end{split}
\label{eq:step3_Fi}
\end{align}

%%% 4
\paragraph{Step 4:} 
As a last step, the remaining free parameters are chosen.
The parameters $K_j$, for $j \in \{n+1,\ldots,N\}$, must be chosen sufficiently large, 
\begin{align}
K_j > K_j^{min}, 
\label{eq:K_j_min}
\end{align}
and the parameters $\epsilon_i$, for $i \in \{1,\ldots,n\}$, sufficiently small, 
\begin{align}
0 < \epsilon_i < \epsilon_i^{max},
\label{eq:eps_i_max}
\end{align}
such that it holds that $\mu_i(x(t)) \in \mathbb{R}_{++}^n$, for all $x \in \mathbb{R}_{++}^N$. 
It was proven in Lemma~2 in \citealp{Schittler2013b} that indeed such $K_j^{min}$ and $\epsilon_i^{max}$ exist.

Thereby interaction functions for all $i=1\ldots N$ are now determined in terms of function classes. The choice of the free parameters within the constraints, as given in the last step, provides degrees of freedom that can be exploited to, for example, fit the dynamics of a model to data.

%%%%%%%%%%%%%%%

%%%%%%%%%%%%%%%
\section{Theorem and proof of multistability equivalence}
\label{sec:proof}

%For the construction method presented in~\citet{Schittler2013b}, the following theorem was stated:
Let us briefly recall the theorem stated in~\citet{Schittler2013b}, regarding the construction method presented therein:

%\vspace{0.2cm}

%%% THEOREM
\begin{theorem}~(\citealp{Schittler2013b})\\
%Given a system~\eqref{eq:sys_f}, 
If~\eqref{eq:S_A} fulfills the Assumptions~1 and~2, %as stated in~\citet{Schittler2013b}, Ass.~1 and Ass.~2, 
then for any system~\eqref{eq:sys_F} as defined in Section~\ref{sec:problem}, 
there exist
$K_j^{min}$, $j = n+1,\ldots,N$ such that every system constructed by the procedure
%presented in~\citet{Schittler2013b} %\textit{steps 1-4} 
given by Steps~1-~4 
with 
$K_j > K_j^{min}$,  $j = n+1,\ldots,N$, % $\epsilon_{ij} < \epsilon_{ij}^{max} \ \forall j \in \{1,\ldots,N\}$,
is multistability-equivalent to system \eqref{eq:sys_f}.
\end{theorem}

The theorem was proven also in~\citet{Schittler2013b} %that this construction yields multistability-equivalent systems as by~(i)-(iii) in Def.~\ref{def:ms-eq}, 
under the assumption that all feedback loops in the interaction functions $a(z)$ can be broken by a single-input single-output (SISO) loopbreaking. 
Such a loopbreaking was defined by~\citet{Waldherr2009}, and it was exploited in the proof for the property~(iii).
However, this loopbreaking approach restricts the proof and thus the applicability of the construction method to a subset of GRN systems.

To overcome this drawback, the idea of a multi-input multi-output (MIMO) loopbreaking based on~\citet{Schittler2013b,Waldherr2009} was introduced by~\citet{Jouini2013}, where also the corresponding generalized proof was sketched.
MIMO systems provide a more general class of systems that allows to cover GRNs even with complex interaction structure.
For the proofs of properties~(i) and~(ii), which are independent of the chosen loopbreaking approach, we refer to~\citet{Schittler2013b}.
Let us now elaborate the generalized proof for property~(iii) in detail.

% LEMMATA?

%%%%%% PROOF
%\begin{proof}

%\vspace{0.2cm}
	%%%%%% (iii)
	%\item[(iii)]  
	\vspace*{2.5ex}
	%\noindent (iii) 
	\paragraph{Proof of~(iii)}
		It remains to be shown that the number of unstable modes
		(eigenvalues with positive real part)
		of the high-dimensional system~\eqref{eq:sys_F}, if constructed by the steps~1-4, is equal to the number of unstable modes 
		(eigenvalues with positive real part)
		of the low-dimensional system~\eqref{eq:sys_f}.
		
		The idea of the proof is the following: 
		At first, for the stability analysis the systems are linearized at steady state.
		Then, we perform a loopbreaking on both systems that supplies inputs and outputs while ensuring that they are no unstable pole-zero cancellations.
		From this, we can now write transfer matrices of these systems. 
		Finally, it is shown that the Nyquist curves of these transfer matrices under mild assumptions are that close to each other that both transfer matrices of the loopbroken systems have the same number of zeros in the right half plane.
		Thus, it can be shown that the transfer matrices of the two closed systems have the same number of poles in the right half plane.
		With unstable pole-zero cancellations ruled out by the initial loopbreaking, we can conclude that the Jacobians of both systems have the same number of eigenvalues in the right half plane.
		The proof was presented first for the special case of systems where a single-input single-output loopbreaking is suffcient~(\citealp{Schittler2013b}). Based on this, the proof was generalized to the class of systems with multi-input multi-output loopbreaking~(\citealp{Jouini2013}), as also presented here.
		
		% First, for stability analysis the systems are linearized at steady state.
		The stability of a system at a steady state can be analyzed via the eigenvalues of its Jacobian, that is of the system linearized at this steady state. 
		Therefore, we consider for system~\eqref{eq:sys_f} the system linearized at a particular steady state $z^{\ast}$,
	  \begin{equation}
		\tilde{\Sigma}_f: 
		\dot{\tilde{z}} = \tilde{f}\tilde{z} = \tilde{a}\tilde{z} - \tilde{d}\tilde{z} ,
		\label{eq:linearized_f}
		\end{equation}
		with the matrices
		\begin{align}
		\begin{split}
		\tilde{a} &:= \frac{\partial a(z)}{\partial z}\Big|_{z = z^{\ast}} \\
		\tilde{d} &:= \frac{\partial d(z)}{\partial z}\Big|_{z = z^{\ast}} ,
		\end{split}
		\label{eq:linearized_f_matrices}
		\end{align} 
		and for system~\eqref{eq:sys_F} the linearized system
	  \begin{equation}
		\tilde{\Sigma}_F: 
		\dot{\tilde{x}} = \tilde{F}\tilde{x} = \tilde{A}\tilde{x} - \tilde{D}\tilde{x} ,
		\label{eq:linearized_F}
		\end{equation}
		with
		\begin{align}
		\begin{split}
		\tilde{A} &:= \frac{\partial A(x)}{\partial x}\Big|_{x = x^{\ast}} \\
		\tilde{D} &:= \frac{\partial D(x)}{\partial x}\Big|_{x = x^{\ast}} .
		\end{split}
		\label{eq:linearized_F_matrices}
		\end{align} 
		%
	%
		
		% Second, we perform a loopbreaking on both systems that supplies inputs and outputs while ensuring that they are no unstable pole-zero cancellations.
		Let us now introduce the loopbreaking of both systems, somewhat similar to the approach in~\citet{Waldherr2009}.
		%The idea is to convert each system into an input-output system, which enables stability analysis by means of its transfer matrix: The poles of the loopbroken system can be analyzed via its Nyquist curve. Then, the direct relation between the eigenvalues of the original system and the transfer matrix of the open loop system allows to conclude about the original system's stability. 
		(For details the reader is referred to~\citet{Waldherr2009}.)
		The aim of the loopbreaking is to obtain a system with inputs and outputs, and at the same time to rule out any unstable pole-zero cancellations.
		The former is required to obtain a transfer matrix, while the latter is important to have every unstable eigenvalue supplied by an unstable pole. 
		Therefore, all loops that may possibly produce such cancellations will be broken:
		The self-loops in the degradation term $\tilde{d}\tilde{z}$ ($\tilde{D}\tilde{x}$, respectively) yield solely negative poles, and thus will not alter the stability. Thus, we only break the loops in $\tilde{a}\tilde{z}$ ($\tilde{A}\tilde{x}$, respectively).
		The obtained systems will be referred to as loopbroken systems, to emphasize that they are not really open-loop systems since the degradation self-loops are still maintained.
		
		We perform such a loop breaking in the low-dimensional system, linearized at a particular steady state, $\tilde{\Sigma}^{(r)}_f$. The corresponding loopbroken system reads%~(\citealp{Jouini2013})
		\begin{align}
		\begin{split}
		\dot{\tilde{z}} &= -\tilde{d}\tilde{z} + \tilde{a}u \\
		y &= I_n \tilde{z} ,
		\end{split}
		\label{eq:openloop_f}
		\end{align}
		with the input $u \in \mathbb{R}_+^n$, output $y \in \mathbb{R}_+^n$, and the linearized dynamics and input matrices as in~\eqref{eq:linearized_f_matrices}.
		The term $\tilde{a}u$ reflects the interaction function of the system~\eqref{eq:linearized_f}, but the argument $\tilde{z}$ is replaced by $u$ to realize the loopbreaking. 
		The original linearized system~\eqref{eq:linearized_f} is obtained by closing the loop via setting the input equal to the output, $u = y$. 
				
		We perform a similar loopbreaking in the high-dimensional system $\Sigma_F$, linearized at a particular steady state $x^{\ast} = h(z^{\ast})$.
		This results in a loopbroken system
		\begin{align}
		\begin{split}
		\dot{\tilde{x}} &= -\tilde{D} \tilde{x} + \tilde{A} U \\
		Y &= I_N \tilde{x} ,	
		\end{split}
		\label{eq:openloop_F}
		\end{align}
		with the input $U \in \mathbb{R}_+^N$, output $Y \in \mathbb{R}_+^N$, and the linearized dynamics and input matrices as in~\eqref{eq:linearized_F_matrices}.	Here, the term $\tilde{A}U$ reflects the linearized interaction function of the original system~\eqref{eq:linearized_F}, but with $\tilde{x}$ replaced by $U$ to realize the loopbreaking. 
		The original linearized system~\eqref{eq:linearized_F} is again obtained by closing the loop via setting the input equal to the output, $U = Y$. 
		
		% From this, we can now write transfer matrices of these systems. 
		We can now write the transfer matrices for the obtained systems. 
		The transfer matrix of the loopbroken system~\eqref{eq:openloop_f} is%~(\citealp{Jouini2013})
		\begin{align}
		\begin{split}
		G_f(\lambda) &= I_n (\lambda I_n + \tilde{d})^{-1} \tilde{a}
		\end{split}
		\end{align}
    and for eigenvalues $\lambda_f$ of the corresponding closed-loop system it holds that
    \begin{align}
    \det \left( I_n - G_f(\lambda_f) \right) = 0
    \end{align}
    where we used~\citet{Waldherr2009}, Lemma~2.3, and the properties of the transfer matrix of an multi-input multi-output system.
    %All loops in $a(z)$ are now broken, and the self-loops in $d(z)$ (corresponding to the linear degradation terms) yield solely negative poles.
    %Thereby, we have ruled out the possibility of unstable pole-zero cancellations -- this is required since otherwise the transfer matrix of the closed-loop system may possibly not supply all positive eigenvalues any more.
    Since now we have ensured that all positive eigenvalues (unstable modes) will be detected as positive poles in the transfer matrix of the closed-loop system, that is, as zeros $\lambda_f$ of $\det(I_n - G_f(\lambda))$, we can exploit $\det(I_n - G_f(\lambda))$ to analyze the unstable modes of the system.
		
		Similarly, the transfer matrix of the second loopbroken system~\eqref{eq:openloop_F} is%~(\citealp{Jouini2013})
		\begin{align}
		\begin{split}
		G_F(\lambda) &= I_N (\lambda I_N + \tilde{D})^{-1} \tilde{A}
		\end{split}
		\label{eq:G_F_openloop}
		\end{align}
		and for eigenvalues $\lambda_F$ of the corresponding closed-loop system it is 
		\begin{align}
		\det \left( I_N - G_F(\lambda_F) \right) = 0 .
		\end{align}

		Let us now set in relationship the winding numbers of the Nyquist curves arising from the two transfer matrices, which allows to deduce the number of unstable poles of the corresponding two closed-loop systems.

		We now let $\varepsilon := \min_{\omega} | \det(I_n - G_f(\j \omega) ) |$ denote the minimum distance of this Nyquist curve from the origin.
		That is, as long as the Nyquist curve $\det(I_N - G_F(\j \omega))$ deviates from the Nyquist curve $\det(I_n - G_f(\j \omega))$ less than $\varepsilon$ (in that it lies within a tube of diameter $\varepsilon$), then both Nyquist curves have the same winding number with respect to the origin.
    Due to the assumption that the Jacobian of $\tilde{\Sigma}_f$ has no eigenvalues on the imaginary axis, we have $\varepsilon > 0$.
    Now,
		for all $j: K_j \rightarrow \infty$ the Nyquist curves approach each other, as has been argued in~\citet{Jouini2013} and will be outlined in the following.

				From~\eqref{eq:G_F_openloop}, it is seen that the rows $(n+1)\ldots N$ of $G_F(\lambda)$ are
				\begin{align}
				\left( G_F(\lambda)|_{i,j} \right)_{\substack{i=(n+1)\ldots N\\ j=1\ldots N}}
				&=
				\diag\left( (\lambda+K_i)^{-1}_{i=(n+1)\ldots N} \right)
				\cdot
				\left( \tilde{A}|_{i,j} \right)_{\substack{i=(n+1)\ldots N\\ j=1\ldots N}} .
				\label{eq:G_F_n+1..N}
				\end{align}
				%consist of entries $../(..+K_j)$.
				For letting the degradation rates of the module genes, $\forall i \in\{n+1,\ldots,N\}:\, K_i = K \rightarrow \infty$, each entry of~\eqref{eq:G_F_n+1..N} becomes
				\begin{align}
				G_F(\lambda)|_{i,j} 
				&=
				(\lambda+K_i)^{-1}
				%\cdot
				\tilde{A}|_{i,j}
				\rightarrow 0,
				\label{eq:G_F_n+1..N_lim}
				\end{align}
				%$/(\lambda+K_i) \rightarrow 0$ 
				by Lemma~1 in~\citet{Schittler2013b}.
		Thus as $K \rightarrow \infty$, the expression for the determinant becomes
		\begin{align}
		\begin{split}
		%\displaystyle
		\det(I_N - G_F(\j \omega)) 
		\rightarrow& 
		\det\big(I_n - (G_F(\j \omega)|_{i,j})_{\substack{i=1\ldots n\\j=1\ldots n}}\big) \cdot \det(I_{N-n}) \\
		&= \det\big(I_n - (G_F(\j \omega)|_{i,j})_{\substack{i=1\ldots n\\j=1\ldots n}}\big) \cdot 1  .
		\end{split}
		\label{eq:det_G_F=det_G_f*1}
		\end{align}
		
		With this it remains to investigate the remaining upper left $(n\times n)$ subblock of $G_F(\lambda)$, corresponding to the transmission of signals between master genes.
		Since $I_N$ and $(\lambda I_N + \tilde{D})^{-1}$ in~\eqref{eq:G_F_openloop} are diagonal matrices, this $(n\times n)$ subblock of $G_F(\lambda)$ only depends on the upper left $(n\times n)$ subblock of the input matrix $\tilde{A}$.  
		Due to the construction of $A(x)$ in step~3, 
		this upper left $(n\times n)$ subblock of the input matrix $\tilde{A}$ in~\eqref{eq:G_F_openloop} is
		\begin{align}
		\begin{split}
		(\tilde{A})_{\substack{i=1\ldots n\\j=1\ldots n}} 
		= 
		\Big( \frac{\partial A_i}{\partial x_j}\Big|_{x^{\ast}} \Big)_{\substack{i=1\ldots n\\j=1\ldots n}} 
		&= 
		\Big( \frac{\partial a_i(\mu_i(x))}{\partial \mu_i(x)} \frac{\partial \mu_i(x)}{\partial x_j} \Big|_{x^{\ast}} \Big)_{\substack{i=1\ldots n\\j=1\ldots n}} 
		.
		\end{split}
		\label{eq:A_1..n__1}
		\end{align}
		In this, $x^{\ast}$ can further be substituted via the map $h:\mathbb{R}^n \rightarrow \mathbb{R}^N$, which was defined in~\citet{Schittler2013b} in the proof of~(ii) to map the steady states of the low-dimensional system to the steady states of the high-dimensional system.
		According to its definition~(\citealp{Schittler2013b}), the map $h$ yields in its first $n$ components, $j\in\{1,\ldots,n\}$:  $x^{\ast}_j = (h(z^{\ast}))_j = z^{\ast}_j$.
		Using this, we can further rewrite~\eqref{eq:A_1..n__1} to
		\begin{align}
		\begin{split}
		\Bigg( \frac{\partial a_i(\mu_i(x))}{\partial \mu_i(x)} \frac{\partial \mu_i(x)}{\partial x_j} \Big|_{x^{\ast}} \Bigg)_{\substack{i=1\ldots n\\j=1\ldots n}} 
		&=
		\Bigg( \frac{\partial a_i(z)}{\partial z} \Big|_{z^{\ast}} \sgn\left( \left| \frac{\partial A(x)}{\partial x_j} \right| \right) \Bigg)_{\substack{i=1\ldots n\\j=1\ldots n}} 
		\\
		&=
		\Big( \frac{\partial a_i(z)}{\partial z_j} \Big|_{z^{\ast}} \Big)_{\substack{i=1\ldots n\\j=1\ldots n}} 
		=
		\frac{\partial a(z)}{\partial z} \Big|_{z^{\ast}}
		= \tilde{a}
		.
		\end{split}
		\label{eq:A_1..n__2}
		\end{align}

		Let us now take this result and the fact that the degradation rates for the master genes, $\forall i\in \{1,\ldots,n\}: K_i = k_i$, into the entries of the $n\times n$ subblock in the transfer matrix, as appearing in~\eqref{eq:det_G_F=det_G_f*1}. Then, we see that % as $K \rightarrow \infty$,
		these entries are
		\begin{align}
		G_F(\j \omega)|_{i,j} = %\rightarrow 
		(\j\omega+k_i)^{-1} \tilde{a}|_{i,j}
		= G_f(\j \omega)|_{i,j}
		.
		\end{align}
		
		Since each entry of the upper $(n\times n)$ subblock of $G_F(\j \omega)|_{i,j}$ equals $G_f(\j \omega)|_{i,j}$, and because of~\eqref{eq:det_G_F=det_G_f*1}, also the determinants approach each other for $K \rightarrow \infty$:
		\begin{align}
		\det(I_N - G_F(\j \omega)) \rightarrow \det(I_n - G_f(\j \omega))
		\label{eq:det_G_F->det_G_f}
		\end{align}
		%
		%For values $K_j < \infty$ the determinants in~\eqref{eq:det_G_F->det_G_f} will deviate but, due to the properties of the limit, there will be a small tube around ...
 
    This in turn implies, similar to the original argument in~\citet{Schittler2013b}, that there exist $K_j^{min}$ such that the Nyquist curve of the high-dimensional system, $\det(I_N-G_F(\mathbf{j}\omega))$, lies within a tube of diameter $\varepsilon = \min_{\omega}\left| \det\left( I_n - G_f(\j\omega) \right) \right|$ around the Nyquist curve of the low-dimensional system
    for sufficiently large $K_j$, and then the winding numbers of these curves around the origin are equal.
                Then, by the argument principle, $\det(I_N-G_F(\mathbf{j}\omega))$ has the same number of zeros in the right half plane as $\det(I_n-G_f(\mathbf{j}\omega))$, and the Jacobians of the two systems $\tilde{\Sigma}_f$ and $\tilde{\Sigma}_F$ have the same number of positive eigenvalues.

With this, it is proven that for sufficiently high $K_j$ the property~(iii) holds, which concludes the proof.
\hspace{0.9\textwidth}
%\begin{flushright}
$\square$
%\end{flushright}
%\end{proof}

%%%%%%%%%%%%%%%
%%%%%%%%%%%%%%%
\section{Example: A GRN in mesenchymal stem cell differentiation}
\label{sec:example}

In this section, we employ the proposed method to examine a GRN %gene regulatory network 
that determines the differentiation of mesenchymal stem cells into the adipogenic, osteogenic, or chondrogenic cell type.
%
% MOTIVATION AND AIM
Mesenchymal stem cells are a type of adult stem cells that are characterized by their potential to differentiate into adipocytes (fat cells), osteoblasts (bone cells), and chondrocytes (cartilage cells)~(\citealp{Baksh2004,Heino2008,Nakashima2003,Ryoo2006}). 

A set of genes and transcription factors have been established as cell type-specific markers for the adipogenic, osteogenic, or chondrogenic lineage, as summarized in Table~\ref{table:MSCgenes}. 
Generally, a core motif model of three key regulator genes that each characterize one of the three considered cell types could be developed based on the dynamical properties of the biological system, for example as in~\citet{Schittler2010} for a similar GRN.
However, such a model neither considers all measured genes, nor can it cope with the higher complexity of the real GRN.
A mathematical model of the detailed GRN, determining mesenchymal stem cell fate, could allow for 
predicting gene expression dynamics under certain differentiation stimuli,
comparing various differentiation protocols, 
or classifying the therapeutic potential of mesenchymal stem cells from individual donors.
Fortunately, the previously presented construction method offers a solution to this discrepancy.

The aim of the example is the following:
\begin{itemize}
\item A low-dimensional GRN model is developed with a parametrization such that it reproduces the observed cell types.
This low-dimensional model is used to study generic properties of the cell differentiation process, such as the stability of cell types in dependence of parameters.
\item A high-dimensional GRN model, obtained via the construction method proposed in this chapter beforehand, is developed that incorporates the set of genes given by the differentiation assays.
\end{itemize}
The obtained high-dimensional model can be fit to the full readout of the experimental data, and can be further exploited to determine for example donor-specific differentiation parameters.
Although the availability of the high-dimensional model offers a whole new field of model application, including parameter estimation, donor classification, and experimental design, these topics are out of the scope of this report.

%\subsection{Characterization of cell types}
%\subsection{Specification of GRNs}
\paragraph{Specification of GRNs}
% SPECIFICATION
%
First, the cell types observed in the mesenchymal stem cell system under study are briefly characterized.
% cell types:
The cell types under consideration are the following three, and are defined according to their expression of type-specific genes as reported in the literature~(\citealp{Baksh2004,Heino2008,Ryoo2006,Nakashima2003,Darlington1998,Tang2004,Drissi2000,Shui2003,Fu2007,Zhou2006}; and others):
\begin{itemize}
\item Adipogenic cell type: This cell type exhibits a high expression at early stages of CEBP$\beta$, followed by CEBP$\alpha$, PPAR$\gamma$, and, at later stages, also LPL. 
\item Osteogenic cell type: This cell type is characterized by early expression of RUNX2, followed by OSX, BGLAP, and SPARC.
\item Chondrogenic cell type: This cell type is characterized by the expression of SOX9. 
\end{itemize}
The cell type-specific genes that have been used in the considered differentiation assays 
have been selected as they have been established in the literature, 
and are also summarized in Table~\ref{table:MSCgenes}.

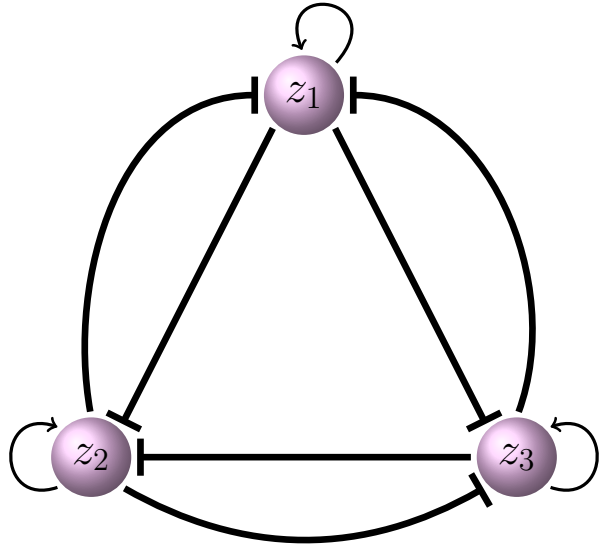
\begin{figure}[!ht]
\begin{center}
\begin{tikzpicture}[font=\Large,scale=0.4,smooth]
    \tikzstyle{shape arrow2}+=[very thick, shorten <=2pt, shorten >=2pt]
		  \tikzstyle{shape arrow3}+=[very thick, shorten <=2pt, shorten >=2pt]
  \tikzstyle{shape arrow4}+=[very thick, shorten <=2pt, shorten >=2pt]
	\tikzstyle{shape arrow5}+=[very thick, shorten <=2pt, shorten >=2pt,line width=2.5pt]
    \tikzstyle{shape mycircle}+=[circle,shading=ball,ball 
color=magenta!20,minimum size=30]
\tikzstyle{shape mycircle2}+=[circle,shading=ball,ball 
color=brown!20,minimum size=24] 
 \begin{scope}      
%CEBPbeta
    \draw (7,12) node[shape mycircle] (x1) {$z_1$}; 
%Sox9
 \draw (14,0) node[shape mycircle] (x3) {$z_3$};
%Runx2
\draw (0,0) node[shape mycircle] (x2) {$z_2$};
 \end{scope}

 %interaktionen
 \begin{scope}
 %Selbstaktivierung
\draw[shape arrow2,->] (x1.north east) .. controls +(2,2.3) and 
+(-1,2.5) .. (x1.north);
 \draw[shape arrow2,->] (x3.south east) .. controls +(2.2,-.75) and 
+(2.2,.75) .. (x3.north east);
\draw[shape arrow2,->] (x2.south west) .. controls +(-2.2,-.75) and 
+(-2.2,.75) .. (x2.north west);
%interaktionen
 \draw[shape arrow5,-|] (x3.west) to (x2.east);%[out=150,in=30]
 \draw[shape arrow5,-|] (x2.south east) to[out=-30,in=210] (x3.south west);%
 \draw[shape arrow5,-|] (x3.north) to [out=70,in=0] (x1.east);
 \draw[shape arrow5,-|] (x1.south east) to (x3.north west);
 \draw[shape arrow5,-|] (x1.south west) to (x2.north east);
 \draw[shape arrow5,|-] (x1.west) to[out=180,in=100] (x2.north);
%zwischen modul master
    \end{scope}		
\end{tikzpicture}
\end{center}
	\caption{Interaction network structure of the low-dimensional GRN, with the state variables representing the master genes as given in Table~\ref{table:MSCgenes}.
	Stump arrows denote negative interactions, sharp arrows denote positive interactions, as given by the interaction sign matrix.}
	\label{fig:MSC_n} 
\end{figure}

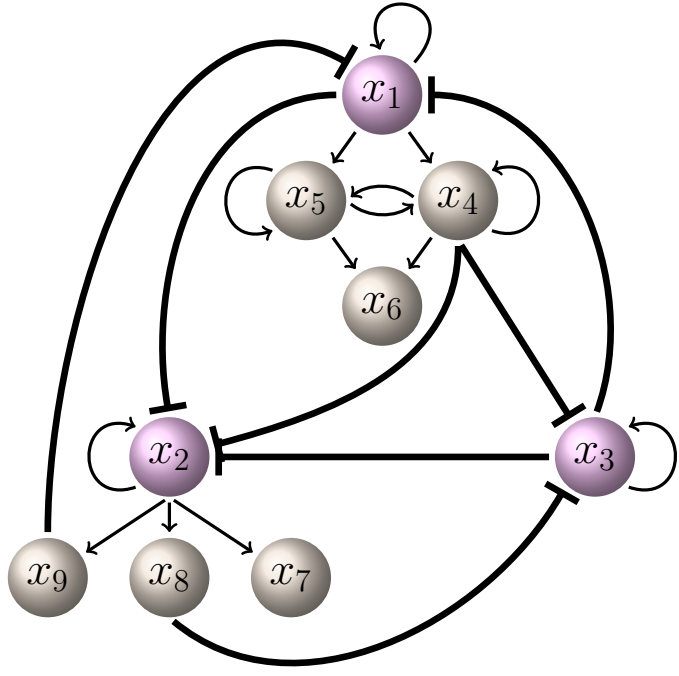
\begin{figure}[!ht]
\begin{center}
\begin{tikzpicture}[font=\Large,scale=0.4,smooth]
    \tikzstyle{shape arrow2}+=[very thick, shorten <=2pt, shorten >=2pt]
		  \tikzstyle{shape arrow3}+=[very thick, shorten <=2pt, shorten >=2pt]
  \tikzstyle{shape arrow4}+=[very thick, shorten <=2pt, shorten >=2pt]
	\tikzstyle{shape arrow5}+=[very thick, shorten <=2pt, shorten >=2pt,line width=2.5pt]
    \tikzstyle{shape mycircle}+=[circle,shading=ball,ball 
color=magenta!20,minimum size=30]
\tikzstyle{shape mycircle2}+=[circle,shading=ball,ball 
color=brown!20,minimum size=24] 
    \begin{scope}  
    
%CEBPbeta
    \draw (7,12) node[shape mycircle] (x1) {$x_1$};
%adipogenes Modul
\draw (9.5,8.5) node[shape mycircle2] (x4) {$x_4$};
 \draw (4.5,8.5) node[shape mycircle2] (x5) {$x_5$};
 \draw  (7,5) node[shape mycircle2] (x6) {$x_6$};
 
%Sox9
 \draw (14,0) node[shape mycircle] (x3) {$x_3$};
		
%Runx2
\draw (0,0) node[shape mycircle] (x2) {$x_2$};
%osteogenes Modul
\draw (0,-4) node[shape mycircle2](x8) {$x_8$};
\draw(-4,-4)  node[shape mycircle2](x9) {$x_9$};
\draw (4,-4) node[shape mycircle2](x7) {$x_7$};
 \end{scope}

 %interaktionen
 \begin{scope}
 %Selbstaktivierung
\draw[shape arrow2,->] (x1.north east) .. controls +(2,2.3) and 
+(-1,2.5) .. (x1.north);
 \draw[shape arrow2,->] (x3.south east) .. controls +(2.2,-.75) and 
+(2.2,.75) .. (x3.north east);
\draw[shape arrow2,->] (x2.south west) .. controls +(-2.2,-.75) and 
+(-2.2,.75) .. (x2.north west);
%zwische master
 \draw[shape arrow5,-|] (x3) to  (x2.east);
 \draw[shape arrow5,-|] (x3.north) to [out=70,in=0] (x1.east);
 \draw[shape arrow5,-|] (x1.west) to[out=180,in=100] (x2.north);
%zwischen modul master
\draw[shape arrow4,->] (x2.south) to (x9);
\draw[shape arrow4,->] (x2.south) to(x8);
\draw[shape arrow4,->] (x2.south) to(x7);
\draw[shape arrow4,->] (x1) to(x4);
\draw[shape arrow5,-|] (x4.south) to(x3);	
\draw[shape arrow5,-|] (x4.south) to[out=-90,in=15](x2);	
\draw[shape arrow4,->] (x1) to(x5);	
\draw[shape arrow3,->] (x4) to(x6);	
\draw[shape arrow3,->] (x5) to(x6);
\draw[shape arrow3,->] (x4.west) to[out=140,in=40] (x5.east);
\draw[shape arrow3,->] (x5.east) to[out=-40,in=-140] (x4.west) ;
\draw[shape arrow5,-|] (x8.south) to[out=-40,in=-120] (x3.south west);
\draw[shape arrow5,-|] (x9.north) to[out=90,in=150] (x1.north west);
%Aktivieerung modul
  \draw[shape arrow2,->] (x4.south east) .. controls +(2.2,-.75) and 
+(2.2,.75) .. (x4.north east);
\draw[shape arrow2,->] (x5.north west) .. controls +(-2.2,.75) and 
+(-2.2,-.75) .. (x5.south west);
    \end{scope}		
\end{tikzpicture}
\end{center}
	\caption{Interaction network structure of the high-dimensional GRN, with the state variables representing all genes (master and module genes) as given in Table~\ref{table:MSCgenes}.
	Master genes are denoted in pink, module genes in grey.
	Stump arrows denote negative interactions, sharp arrows denote positive interactions, as given by the interaction sign matrix.}
	\label{fig:MSC_N} 
\end{figure}

The genes CEBP$\beta$, RUNX2, and SOX9 are established as master regulators in the cell differentiation of mesenchymal stem cells into adipogenic, osteogenic, and chondrogenic lineage, respectively, and are known to regulate also the expression of the other genes that were selected (Cf. references above).
Based on literature research again, two interaction networks can be drawn: The interaction network of a low-dimensional GRN between the $n=3$ master genes as depicted in Figure~\ref{fig:MSC_n}, and the interaction network of a higher-dimensional GRN between all $N=9$ genes, as depicted in Figure~\ref{fig:MSC_N}.
Therein, stump arrows denote negative interactions, sharp arrows denote positive interactions, as defined by interaction sign matrix entries~\eqref{eq:S_a}.

\begin{table}[!t]
\centering
\begin{tabular}{|c|c|c|c|}
\hline
Gene  & full or alternative name  &   cell type-specificity   &  state variable \\
\hline
\hline
CEBP$\beta$  & CCAAT-enhancer-binding protein $\alpha$ & adipogenic  & $z_1$, $x_1$  \\
RUNX2        & Runt-related transcription factor 2 & osteogenic  & $z_2$, $x_2$  \\
SOX9    		 & Sry-related HMG box & chondrogenic & $z_3$, $x_3$ \\
\hline
PPAR$\gamma$ & Peroxisome proliferator-activated receptor $\gamma$ & adipogenic  & $x_4$ \\
CEBP$\alpha$ & CCAAT-enhancer-binding protein $\beta$ & adipogenic  &  $x_5$ \\
LPL          & Lipoprotein lipase & adipogenic  & $x_6$ \\
\hline
BGLAP        & Osteocalcin %Bone gamma-carboxyglutamate protein %/ Osteocalcin 
		& osteogenic  & $x_7$ \\
OSX          & Osterix & osteogenic & $x_8$ \\
SPARC        & Osteonectin %Secreted protein acidic and rich in cysteine %/ Osteonectin 
    & osteogenic  & $x_9$ \\
\hline
\end{tabular}
\caption{Genes that are established as adipogenic, osteogenic, or chondrogenic markers respectively, of which the mRNA products have been measured in the considered differentiation assays of mesenchymal stem cells. 
The first three genes denote ``key regulators'' or master genes, whereas the remaining genes are additionally measured.} 
\label{table:MSCgenes}
\end{table}

Let the corresponding dynamical GRN model of the low-dimensional GRN be represented by the following system of ODEs, with the state variables $z \in \mathbb{R}^3$ encoding the expression levels of master genes as summarized in Table~\ref{table:MSCgenes}:
\begin{align}
\begin{split}
\dot{z}_1 &= \frac{0.2 z_1^2 + 0.5 + u_A}{10m + 0.1 z_1^2 + 0.5 z_2^2 + 0.5 z_3^2} - 0.1 z_1
\\
\dot{z}_2 &= \frac{0.1 z_2^2 + 1 + u_O}{1m + 0.1 z_2^2 + 0.5 z_1^2 + 0.1 z_3^2} - 0.1 z_2
\\
\dot{z}_3 &= \frac{0.1 z_3^2 + 1 + u_C}{1m + 0.1 z_3^2 + 0.5 z_1^2 + 0.1 z_2^2} - 0.1 z_3
\end{split}
\label{eq:MSC_n}
\end{align}

\begin{figure}[!ht]
\centering
\includegraphics[width=0.65\textwidth]{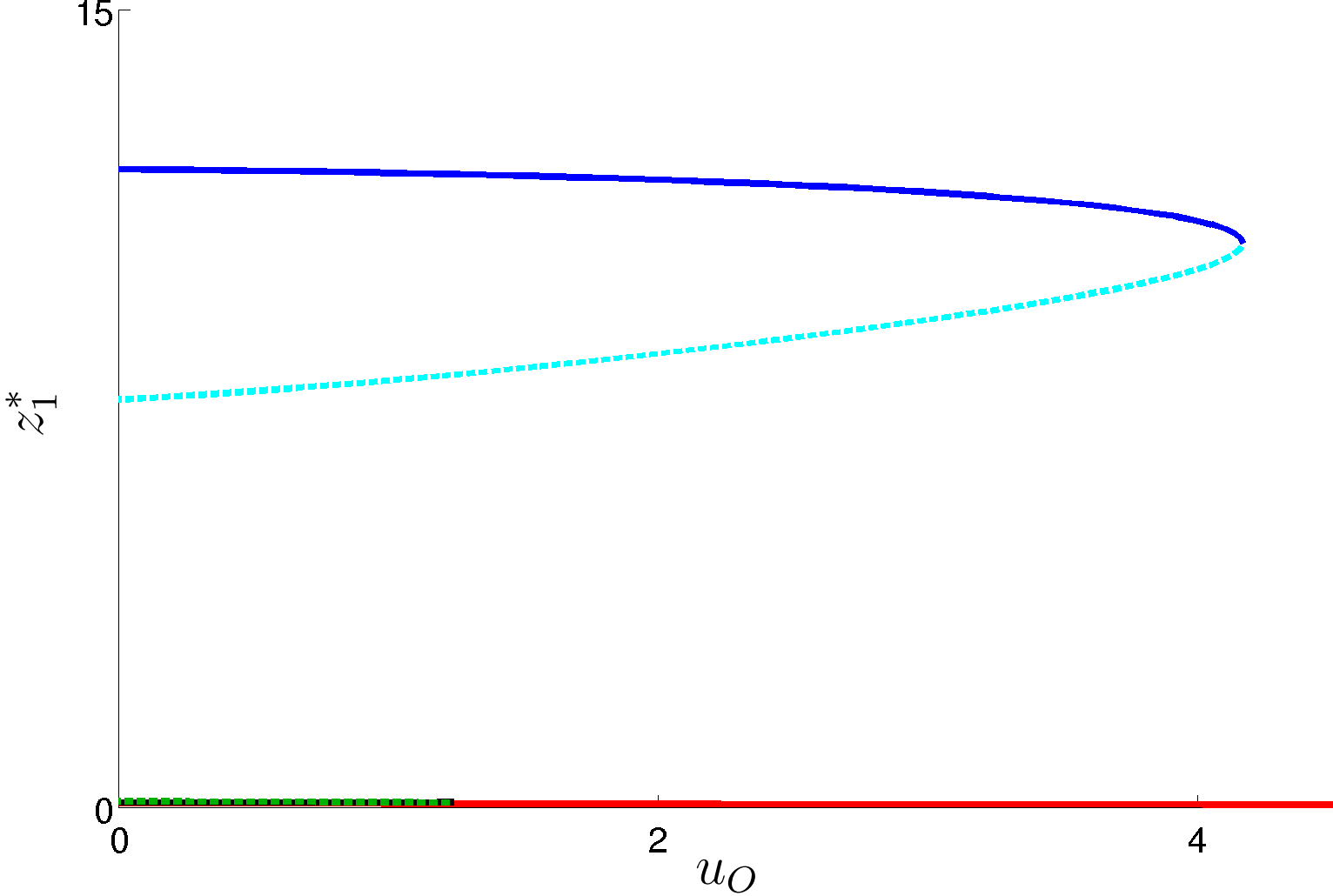}\\
\vspace{-0.65cm}
\includegraphics[width=0.65\textwidth]{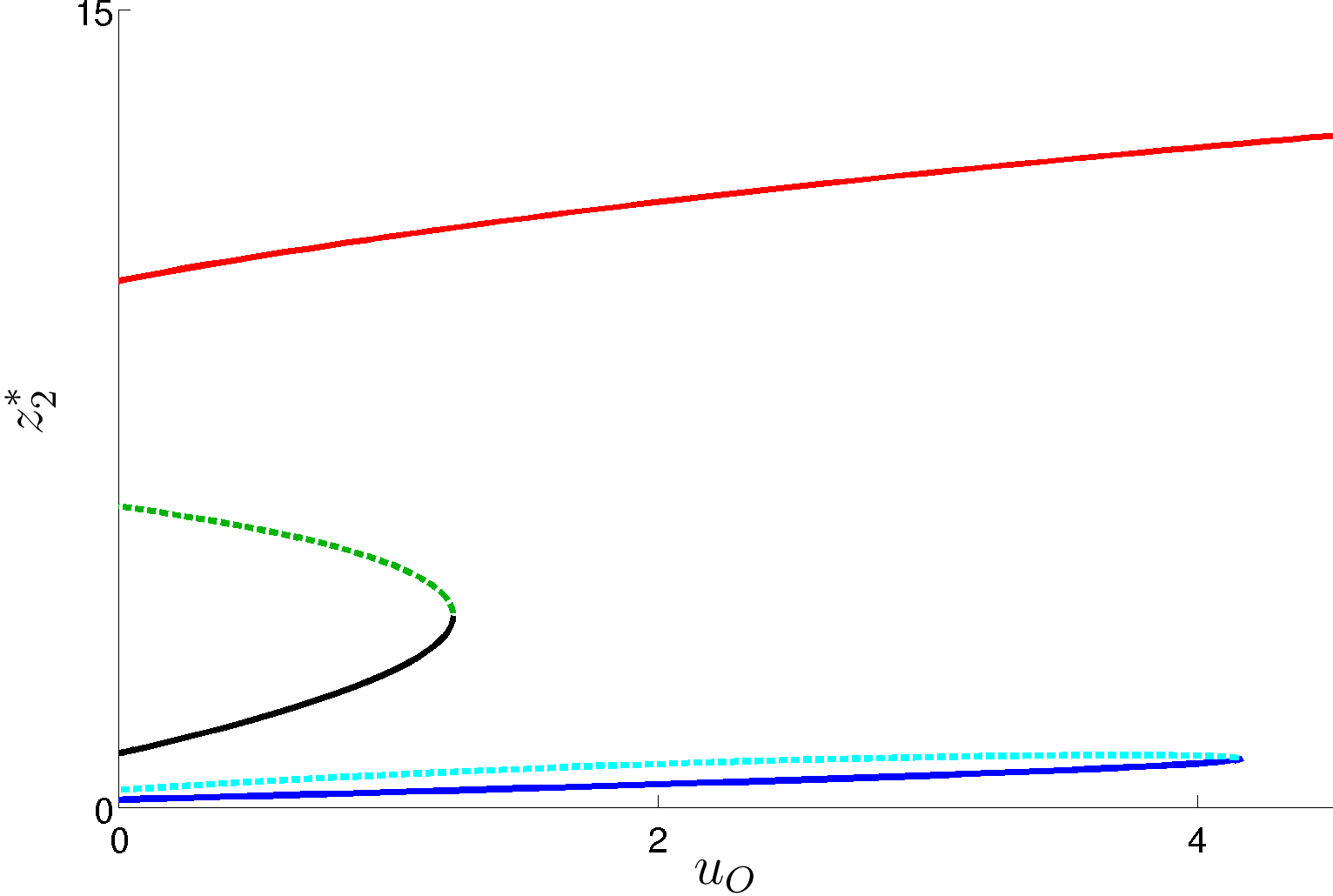}\\
\vspace{-0.65cm}
\includegraphics[width=0.65\textwidth]{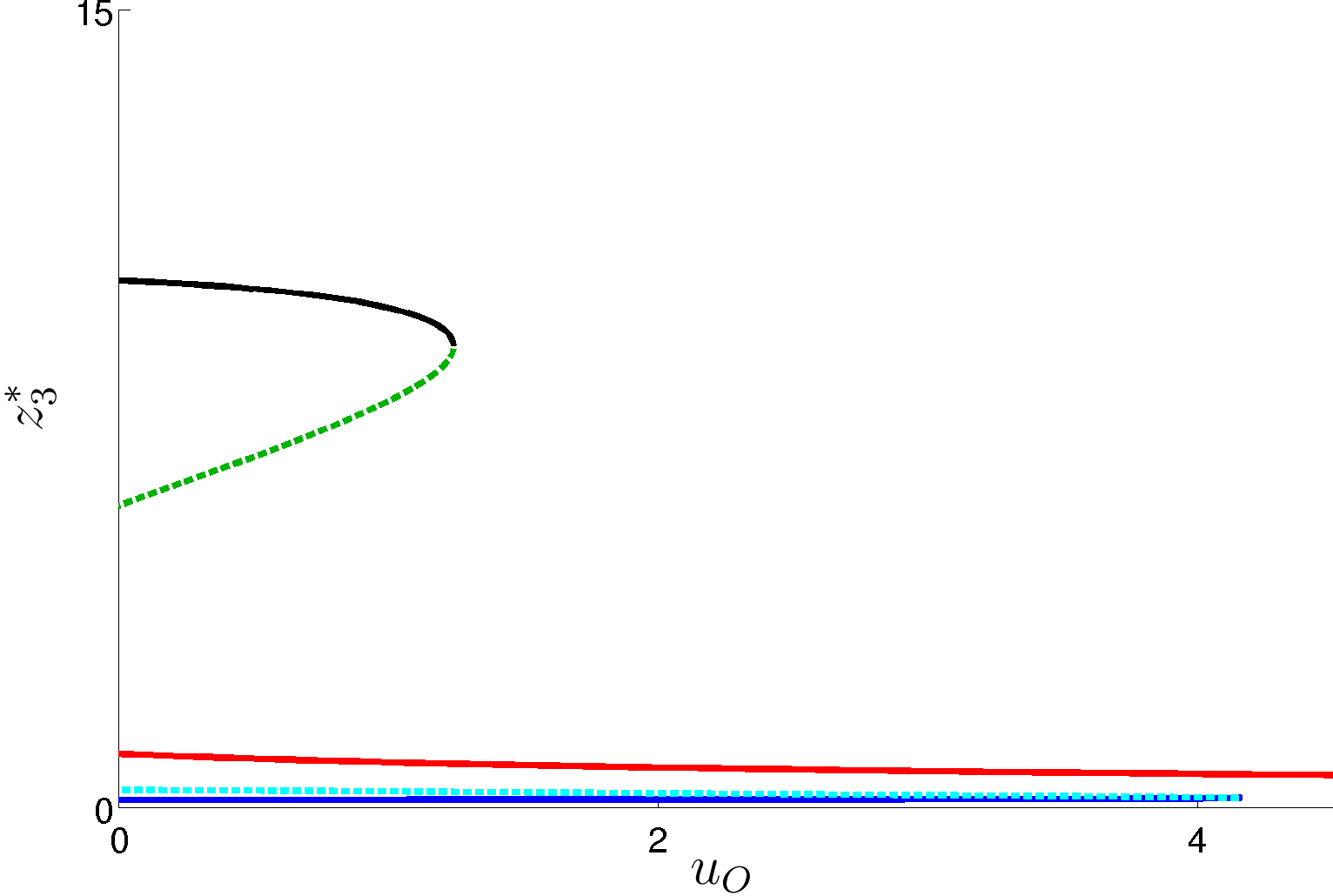}
\vspace{-0.4cm}
\caption{Bifurcation analysis of the low-dimensional system:
If no osteogenic stimulus is applied ($u_O=0$), there are three stable steady states (red, black, blue) and two unstable steady states (green dotted, cyan dotted). The three stable steady states correspond to three cell types: An adipogenic cell type (blue; high $x_1$, low $x_2, x_3$), an osteogenic cell type (red; high $x_2$, low $x_1, x_3$), and a chondrogenic cell type (black; high $x_3$, low $x_1, x_2$).
Upon increasing the osteogenic stimulus $u_O$, first the chondrogenic cell type vanishes, then also the adipogenic cell type vanishes, each by a saddle-node bifurcation. Thus, for high values of osteogenic stimulus $u_O>4.2$, only the osteogenic cell type remains as a stable steady state, and thus the system will converge towards the osteogenic cell type.
If the osteogenic stimulus is withdrawn, this stable cell type is maintained.
}
\label{fig:MSC_bifuana}
\end{figure}

From the observations of the biological system, the following system properties should be captured by the model:
%The parameters have been chosen such that: 
\begin{itemize}
  \item[(S1)] There are four free parameters, corresponding to: suppression from stem cell maintenance factors, $m$, adipogenic stimulus, $u_A$, osteogenic stimulus, $u_O$, chondrogenic stimulus, $u_C$. % --> also in Table?!
	\item[(S2)] If the system is unstimulated ($u_O=u_A=u_C = 0$) and the stem cell maintenance is low ($m=1$), the model exhibits three stable steady states, corresponding to the adipogenic, osteogenic, chondrogenic cell type with gene expression levels as described above.
	\item[(S3)] If an adipogenic (or, osteogenic, chondrogenic) stimulus with sufficiently high value is applied, $u_A > u_A^{crit}$ ($u_O > u_O^{crit}$, $u_C > u_C^{crit}$, respectively), only one stable steady state remains which corresponds to the adipogenic (osteogenic, chondrogenic) cell type.
	\item[(S4)] If the suppression from the stem cell maintenance factor is kept at a sufficiently high value, $m > m^{crit}$, only one stable steady state remains which corresponds to the mesenchymal stem cell type with low gene expression levels of all type-specific genes.
\end{itemize}
These properties are fulfilled for the chosen parameter values, as was ensured via numerical solution of the system of equations $\dot{z}=0$ and bifurcation analysis, for the parameters $u_A, u_O, u_C, m$, respectively. 
Exemplarily, bifurcation analysis results are shown for the parameter $u_O$, the osteogenic stimulus, in Figure~\ref{fig:MSC_bifuana}, and discussed later.

%\subsection{Construction of high-dimensional GRN}
\paragraph{Construction of high-dimensional GRN}

Given the interaction network of the high-dimensional GRN~\ref{fig:MSC_N}, and the dynamics of the low-dimensional GRN~\eqref{eq:MSC_n}, the construction procedure presented in Section~\ref{sec:construction} and in~\citet{Schittler2013b} can now be conducted.

The obtained steady state gain parameters are
\begin{align}
\begin{split}
\gamma_{41} &= \frac{K_5}{K_4 K_5 - K_4 - K_5}, 
\gamma_{51} = \frac{K_4}{K_4 K_5 - K_4 - K_5}, 
\gamma_{61} = \frac{K_4 + K_5}{K_6 (K_4 K_5 - K_4 - K_5)}, 
\\
\gamma_{72} &= \frac{1}{K_7}, \hspace{2.65cm}
\gamma_{82} = \frac{1}{K_8},  \hspace{2.65cm}
\gamma_{92} = \frac{1}{K_9}.
\end{split}
\label{eq:MSC_gammas}
\end{align}

This yields the following system of ODEs giving the dynamics of the high-dimensional GRN, with state variables $x \in \mathbb{R}^9$:
\begin{align}
\begin{split}
\dot{x}_1 &= \frac{0.2 x_1^2 + 0.5 + u_A}{10m + 0.1 x_1^2 + 0.5 (\gamma^{-1}_{92}x_9)^2 + 0.5 x_3^2} - 0.1 x_1
\\
\dot{x}_2 &= \frac{0.1 x_2^2 + 1 + u_O}{1m + 0.1 x_2^2 + 0.5 (\frac{x_1+\gamma^{-1}_{41}x_4}{2})^2 + 0.1 x_3^2} - 0.1 x_2
\\
\dot{x}_3 &= \frac{0.1 x_3^2 + 1 + u_C}{1m + 0.1 x_3^2 + 0.5 (\gamma^{-1}_{41}x_4)^2 + 0.1 (\gamma^{-1}_{82}x_8)^2} - 0.1 x_3
\\
\dot{x}_4 &= x_1 + x_4 + x_5 - K_4 x_4
\\
\dot{x}_5 &= x_1 + x_4 + x_5 - K_5 x_5
\\
\dot{x}_6 &= x_4 + x_5 - K_6 x_6
\\
\dot{x}_7 &= x_2 - K_7 x_7
\\
\dot{x}_8 &= x_2 - K_8 x_8
\\
\dot{x}_9 &= x_2 - K_9 x_9
\end{split}
\label{eq:MSC_N}
\end{align}
with $\gamma_{ij}$ as in~\eqref{eq:MSC_gammas}, and $K_j$ as free parameters to be chosen.

For the remainder, we choose $K_4 = K_5 = 3, K_6 = K_7 = K_8 = K_9 = 1$. %, if not mentioned other.
Furthermore, the default values for the other parameters are $m=1$, $u_A=u_O=u_C=0$, corresponding to low levels of the stem cell maintenance factor, and no application of cell type specific stimuli.

%\subsection{Multistability and bifurcation analysis}
\paragraph{Multistability and bifurcation analysis}

% steady states + eigenvalues:
The steady states of the low-dimensional system were obtained by numerically solving the system of equations $\dot{z}=0$.
They are listed, along with their eigenvalues, in Table~\ref{tab:SSs_EWs_n}.
The first three steady states are, according to their eigenvalues, stable, and correspond to the adipogenic, osteogenic, and chondrogenic cell type, as defined by their gene expression levels above.
The steady states of the high-dimensional system, along with their eigenvalues, are depicted in Table~\ref{tab:SSs_EWs_N}.
From there it can be seen that the high-dimensional system indeed shares the same multistability properties as the low-dimensional system.
While the steady states equal in their first $n=3$ components, %(as also reflected by the map $h$...),
the eigenvalues may slightly differ in their exact values but share the same sign and are, for the chosen parameters, close to the eigenvalues of the low-dimensional system.

\begin{table}[!tb]
\centering
\begin{tabular}{| l | ccccc |}
\hline
(a)					& $z^{\ast,1}$	& $z^{\ast,2}$	& $z^{\ast,3}$	& $z^{\ast,4}$	& $z^{\ast,5}$ \\
\hline
$z_1^{\ast(r)}$	&   12.00 &  0.08  &  0.08  &  7.67  &  0.12 \\
$z_2^{\ast(r)}$	&    0.14 &  9.90  &  1.01  &  0.33  &  5.67 \\
$z_3^{\ast(r)}$	&    0.14 &  1.01  &  9.90  &  0.33  &  5.67 \\
\hline
\end{tabular}
\\
$\;$
\\
\begin{tabular}{| l | ccccc |}
\hline
(b)					& $\lambda_f^{1}$	& $\lambda_f^{2}$	& $\lambda_f^{3}$	& $\lambda_f^{4}$	& $\lambda_f^{5}$ \\
\hline
$\lambda_1^{(r)}$ &  -0.02 & -0.11 & -0.10 & +0.02 & -0.10 \\
$\lambda_2^{(r)}$ &  -0.10 & -0.10 & -0.11 & -0.10 & -0.12 \\
$\lambda_3^{(r)}$ &  -0.10 & -0.07 & -0.07 & -0.10 & +0.05 \\
\hline
\end{tabular}
\caption{Steady states (a) and their eigenvalues (b) of the low-dimensional system.}
\label{tab:SSs_EWs_n}
\end{table}

\begin{table}[!tb]
\centering
\begin{tabular}{| l | ccccc |}
\hline
(a)					& $x^{\ast,1}$	& $x^{\ast,2}$	& $x^{\ast,3}$	& $x^{\ast,4}$	& $x^{\ast,5}$ \\
\hline
$x_1^{\ast(r)}$	&   12.00 &  0.08  &  0.08  &  7.67  &  0.12 \\
$x_2^{\ast(r)}$	&    0.14 &  9.90  &  1.01  &  0.33  &  5.67 \\
$x_3^{\ast(r)}$	&    0.14 &  1.01  &  9.90  &  0.33  &  5.67 \\
$x_4^{\ast(r)}$	&   12.00 &  0.08  &  0.08  &  7.67  &  0.12 \\
$x_5^{\ast(r)}$	&   12.00 &  0.08  &  0.08  &  7.67  &  0.12 \\
$x_6^{\ast(r)}$	&   24.00 &  0.17  &  0.17  &  15.35 &  0.24 \\
$x_7^{\ast(r)}$	&    0.14 &  9.90  &  1.01  &  0.33  &  5.67 \\
$x_8^{\ast(r)}$	&    0.14 &  9.90  &  1.01  &  0.33  &  5.67 \\
$x_9^{\ast(r)}$	&    0.14 &  9.90  &  1.01  &  0.33  &  5.67 \\
\hline
\end{tabular}
\\
$\;$
\\
\begin{tabular}{| l | ccccc |}
\hline
(b)					& $\lambda_F^{1}$	& $\lambda_F^{2}$	& $\lambda_F^{3}$	& $\lambda_F^{4}$	& $\lambda_F^{5}$ \\
\hline
$\lambda_1^{(r)}$ &  -0.02 & -0.11 & -0.10 & +0.02 & -0.10 \\
$\lambda_2^{(r)}$ &  -0.10 & -0.10 & -0.11 & -0.10 & -0.13 \\
$\lambda_3^{(r)}$ &  -0.10 & -0.07 & -0.07 & -0.10 & +0.05 \\
$\lambda_4^{(r)}$ &  -1.00 & -1.00 & -1.00 & -1.00 & -1.00 \\
$\lambda_5^{(r)}$ &  -1.00 & -1.00 & -1.00 & -1.00 & -1.00 \\
$\lambda_6^{(r)}$ &  -1.00 & -1.00 & -1.00 & -1.00 & -1.00 \\
$\lambda_7^{(r)}$ &  -1.00 & -1.00 & -1.00 & -1.01 & -1.00 \\
$\lambda_8^{(r)}$ &  -1.00 & -1.00 & -1.00 & -1.00 & -1.00 \\
$\lambda_9^{(r)}$ &  -3.00 & -3.00 & -3.00 & -3.00 & -3.00 \\
\hline
\end{tabular}
\caption{Steady states (a) and their eigenvalues (b) of the high-dimensional system.}
\label{tab:SSs_EWs_N}
\end{table}

In order to investigate the changes in steady states upon changes in the parameters, we conducted bifurcation analysis of the low-dimensional system (via the software package CL\_MatCont, \citet{Dhooge2003}).
The bifurcation diagrams in Figure~\ref{fig:MSC_bifuana} depict the effect of an osteogenic stimulus $u_O$: As long as no osteogenic stimulus is applied ($u_O=0$), there are three stable steady states, and two unstable steady states, as also predicted by the steady state analysis above.
For a sufficiently high osteogenic stimulus ($u_O>u_O^{crit} \approx 4.2$), only the osteogenic cell type remains as a stable steady state.
Thus, the system will converge towards the osteogenic cell type, which is maintained when the stimulus is withdrawn.
Analogous results hold for the effects of an adipogenic stimulus, $u_A$, and a chondrogenic stimulus, $u_C$ (results not shown).

These results are in accordance with the system properties~(S1)-(S3) that were required to be captured by the model.
In addition, the effect of a stem cell maintenance factor $m$ can be investigated similarly. 
Bifurcation analysis reveals that, for a low stem cell maintenance factor $m=1$, there are the three stable cell types as mentioned, and seen in Figure~\ref{fig:MSC_bifuana}.
For high levels of the stem cell maintenance factor $m>m^{crit} \approx 4.5$, however, only one stable steady state remains with low levels in all three cell type specific genes, corresponding to a stem cell state (results not shown).

The bifurcation analysis conducted on the low-dimensional system exemplifies how the dependency of multistability properties on parameters can be investigated in the core motif model.
As a result of the construction procedure, ensuring multistability equivalence of the high-dimensional system, we know that also the high-dimensional system will have the same multistability properties for the discussed high values of the parameters $u_O$, $m$, and also $u_A$ and $u_C$.
It has to be emphasized that, at a saddle-node or pitchfork bifurcation point, there will be an imaginary eigenvalue, thus the multistability equivalence will not hold at exactly this point.
However, for parameter values that are not in the critical region, the multistability equivalence is fulfilled, and we can conclude about the multistability properties of the high-dimensional system from the multistability properties of the low-dimensional system.
In this way, we have circumvented the need for an exhaustive bifurcation analysis of the high-dimensional system.

%\begin{landscape}
%\end{landscape}

%\subsection{Dynamical simulations}
\paragraph{Dynamical simulations}

\begin{figure}[!ht]
\centering
\includegraphics[width=0.72\textwidth]{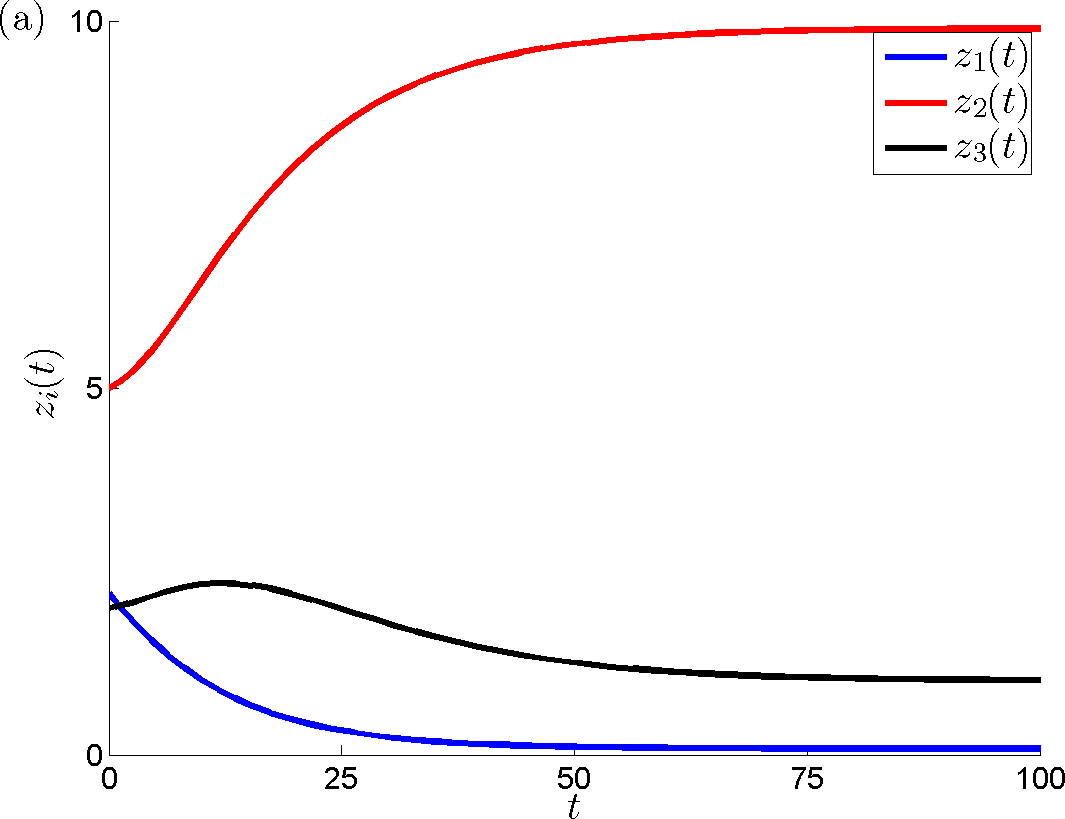}\\
%\vspace{-0.65cm}
\includegraphics[width=0.72\textwidth]{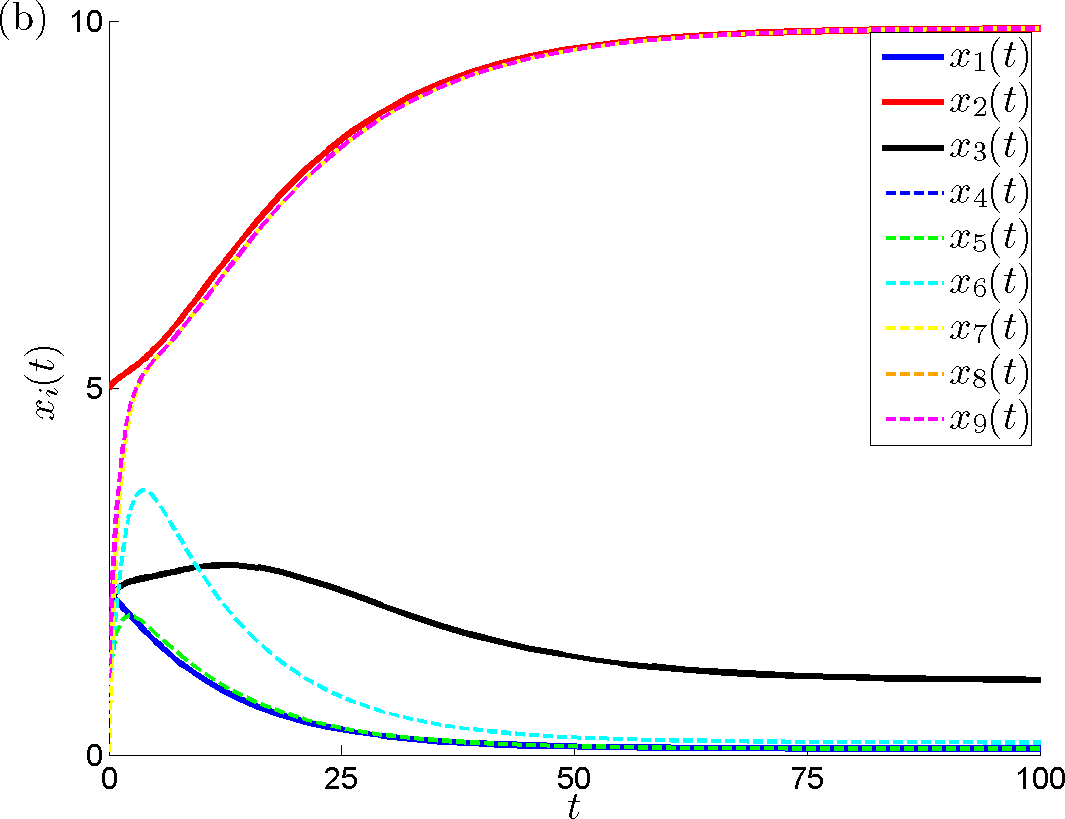}
\vspace{-0.2cm}
\caption{Simulations of dynamics of the low-dimensional system (a) and of the high-dimensional system (b).
Parameters and initial values as in the text.
The dynamics of the master genes in the high-dimensional system are very similar to the dynamics in the low-dimensional system.
}
\label{fig:MSC_dynamics}
\end{figure}

Exemplary simulations of dynamics are shown in Figure~\ref{fig:MSC_dynamics}, for initial values $z(0) = [2.2,5,2]^T$, and $x(0)=[2.2,5,2,1,1,1,0,1,1]^T$, respectively.
As can be seen, the dynamics of the master genes in the high-dimensional system are quite similar to the dynamics in the low-dimensional system, as has been observed also in simulations for other initial values.
This similarity between the dynamics could already be expected from the high similarity in the eigenvalues of the systems, recalling Tables~\ref{tab:SSs_EWs_n} and~\ref{tab:SSs_EWs_N}.
%The higher the degradation rates of the module genes $K_j$, the closer the module gene dynamics come to a quasi steady state.
For choosing very high values for the degradation rate parameters of the module genes, $K_j,\, j\in \{4,\ldots,9\}$, the dynamics of the master genes approach the dynamics of the low-dimensional system even more (results not shown).
% approach infinity, $K_j \rightarrow \infty,\, j\in \{4,\ldots,9\}$
For such parameter values it was furthermore shown in a similar example by~\citet{Jouini2013} that for each steady state of the high-dimensional system, the Nyquist curve approaches the Nyquist curve of the corresponding steady state of the low-dimensional system.
These findings on the dynamics and the Nyquist curves both illustrate the property exploited in the proof of~(iii), namely that by letting $K_j  \rightarrow\infty,\, j\in \{n+1,\ldots,N\}$, the Nyquist curves of the two systems can be made arbitrarily close.

%%%%%%%%%%%%%%%
%%%%%%%%%%%%%%%
\section{Summary and outlook}
\label{sec:summ}

In this report, we presented the generalization of a construction method for multistability-equivalent GRNs of different dimensionality.
Given the dynamics of a low-dimensional GRN, and the interaction structure of a high-dimensional GRN which represents a more detailed expansion of the former, 
the proposed method serves to construct for the high-dimensional GRN a dynamical model with the same multistability properties, in terms of steady states and their stability.
Here, this method was generalized to GRN systems with arbitrarily many internal feedback loops, by using a loopbreaking approach that leads to a multi-input multi-output (MIMO) system.
With this generalization, our method becomes applicable to a broad class of GRN models which have to fulfill only some mild technical assumptions.

By studying an example of a GRN in mesenchymal stem cell differentiation, we demonstrated the potential and value of our method.
A low-dimensional ``core motif'' GRN served to derive a model which meets the properties as observed in the biological system, and to investigate the effects of parameter changes as induced for example by differentiation stimuli.
A high-dimensional GRN model was then derived via the proposed method and can be fit to the gene expression readout of differentiation experiments.
With this, we showed how the results obtained from low-dimensional core motif models can be transferred to more realistic and detailed high-dimensional models of GRNs.

As the construction method proposed here represents sufficient, but not necessary conditions for multistability equivalence, alternative construction methods may be developed. 
These may well exploit additional degrees of freedom that were not pursued in this work. 
For example, the linear interaction functions for module genes may be replaced by more general nonlinear functions, or at least they may be formulated more generally with interaction parameters $\alpha_{ij} \in \mathbb{R}$ replacing $S_A|_{i,j} \in\{-1,0,+1\}$.
Other topics for future work could be whether and under which circumstances equivalences may also be reasoned for bifurcations, or/and vector fields.
For example, related questions, regarding saddle-node bifurcations in systems of different dimensionality but in the context of model reduction, have been considered by~\citet{ChiangFek1993}.
The here investigated example provides evidence that bifurcations and dynamics may, at least to some extent, have similar properties between multistability-equivalent systems.

Summing up, we have introduced the concept of multistability equivalence between GRN systems of different dimensionality, proposed a construction method which we have proven to yield a multistability-equivalent system, and demonstrated its value by investigating an example GRN. 
Our method contributes to overcoming the gap between modeling approaches on distinct levels of detail.
It opens up new possibilities of integrating the results from multistability analysis into the development and parametrization of realistic GRN models.

%%%%%%%%%%%%%%%
% References
\bibliographystyle{plainnat}
\bibliography{refs_schittler2013_TechRep}

\end{document}